\documentclass{article} 

\usepackage{amsmath,amsthm}     
\usepackage{graphicx}     
\usepackage{hyperref} 
\usepackage{url}
\usepackage{amsfonts} 

\usepackage{ amssymb } 

\usepackage{multicol}

\usepackage{tikz} 
\usepackage{pgfplots} 
\usetikzlibrary{calc} 

\newcommand{\polygon}[2]{%
  let \n{len} = {2*#2*tan(360/(2*#1))} in
 ++(0,-#2) ++(\n{len}/2,0) \foreach \x in {1,...,#1} { -- ++(\x*360/#1:\n{len})}}

\theoremstyle{theorem}
\newtheorem{theorem}{Theorem}

\theoremstyle{definition}
\newtheorem*{definition}{Definition}
\newtheorem*{remark}{Remark}

\allowdisplaybreaks

\makeatletter
\@addtoreset{footnote}{page}
\makeatother

\begin{document}

\title{Approximating Mathematical Constants using Minecraft}



\begin{center}
\huge Approximating Mathematical Constants using Minecraft \\

\begin{multicols}{2}

\Large Molly Lynch\\               
\scriptsize Hollins University\\    
lynchme2@hollins.edu \\

\columnbreak

\Large Michael Weselcouch\\               
\scriptsize Roanoke College\\    
weselcouch@roanoke.edu \\

\end{multicols}


\end{center}

\noindent Imagine that you're trapped in the video game Minecraft and that the only way you can return to the real world is if you successfully approximate the values of various mathematical constants using the mechanics of the game.  Would you be able to escape or would you be trapped in Minecraft forever? 

Minecraft is a sandbox video game that revolves around picking up and placing blocks. These blocks are arranged in a three dimensional grid, but the players can move freely around the world. Within the game, players can find raw materials to craft tools and items as well as build structures. It is the best selling video game of all time with over 238 million copies sold. Due to the versatility of the game, Minecraft has long been used in educational settings. In 2011, Minecraft released Minecraft Education Edition and has a number of lessons developed in subject areas including science, math, and history and culture \cite{minecraftedu}. However, most of what has been done with Minecraft for education has been aimed at K-12 classrooms. Up until now, very little has been done to use Minecraft as a tool to explain and explore college level mathematics. 

In this article we will use Minecraft to experimentally approximate the values of four different mathematical constants. Approximations of some of these numbers date back as early as the Indian Shulva Sutras from the Vedic period ( c. 1500 – c. 500 BCE) \cite{shulvasutras}. In these, the square root of two was approximated for use in the creation of an altar for ritual sacrifice as well as approximations for $\pi$ in work constructing a circle with the same area as a given square. 

The mathematical constants that we will approximate are $\sqrt{2}, \pi$, Euler's number $e$, and Ap\'{e}ry's constant $\zeta(3)$. We will begin each section with a brief history of the number being approximated and describe where it appears in mathematics. We then explain how we used Minecraft mechanics to approximate the constant. At the end of each section, we provide some ideas for how to apply our techniques to the approximation of other mathematical constants in Minecraft or elsewhere.  This article is a proof of concept that Minecraft can be used in higher education.

We tried to not limit ourselves to just one area of mathematics when making the approximations.  In this article, we will apply techniques from geometry to approximate $\sqrt{2}$, techniques from calculus to approximate $\pi$, techniques from combinatorics to approximate $e$, and techniques from number theory to approximate Ap\'{e}ry's constant. Any necessary background information on the mathematics used in our approximations will be given at the beginning of the corresponding section. 

Although Minecraft is the tool we used to run the experiments, we will assume the reader has no previous knowledge of the game. In the next section we provide some details about the Minecraft mechanics that are used in the approximations.
We should note that the goal of this article is not to have the most accurate approximations possible, the goal is to inspire people to have fun while learning about various mathematical topics.  We hope you learn something new in this article and feel inspired to try some of these techniques on your own.

\section{Minecraft's mechanics}

Before we can use Minecraft to approximate these mathematical constants, we first provide background information on some basic Minecraft mechanics and items.  If you are already well-versed in Minecraft or are just interested in the math, feel free to skip this section.
\subsection{The hopper}
\begin{figure}[h]
\begin{center}
\includegraphics[height=5cm,width=4.17cm]{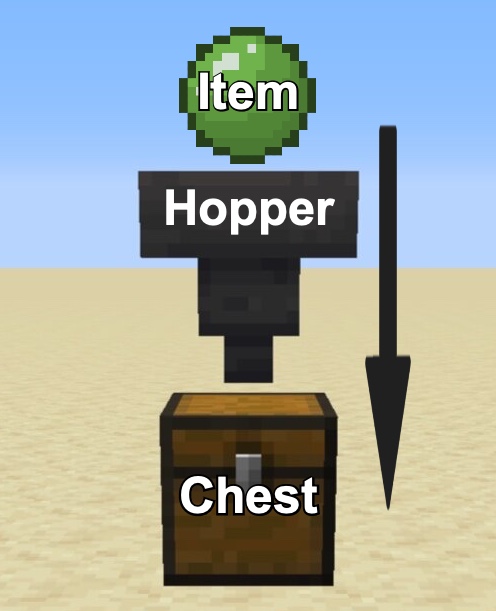}
\caption{A hopper collect items above it and put them in chests.}\label{Hopper}
\end{center}
\end{figure}

A \textit{hopper} is a block in Minecraft that can be used to manage items.  In particular, a hopper collects any item that is directly above it (see Figure \ref{Hopper}). This means it can be used to keep track of location of a dropped item.  If a mob (player, animal, or monster) is killed while standing on a hopper, then the hopper will collect the items dropped by the mob.  Therefore hoppers can be used to record the location where a mob was killed.

In addition to collecting items, another property that hoppers have is that they release items with a constant rate of $2.5$ items per second.  The hopper's ability to release items can be toggled on and off.  Since hoppers release items at a constant rate, they can be used as timers.  For example, if a hopper releases 25 items, then we know that the amount of time that the hopper was releasing items is between $10$ seconds and $10.4$ seconds.  We should note that there are ways to make more accurate timers in Minecraft, but a hopper timer will be accurate enough for our experiments.

\subsection{The dropper}
\begin{figure}[h]
\begin{center}
\includegraphics[height=5cm,width=5cm]{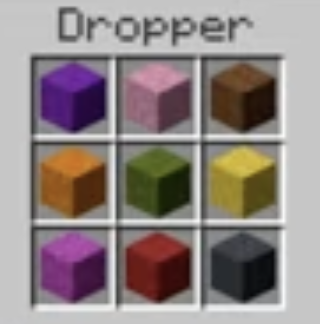}
\caption{A dropper can contain up to nine different blocks.}\label{Dropper}
\end{center}
\end{figure}

A \textit{dropper} is a block that can eject an item.  Droppers have the ability to hold up to $9$ different items at a time (see Figure \ref{Dropper}).  When activated, the dropper randomly selects one of its items to eject.  As a result, droppers can be used as randomizers.  For example, if the dropper has 5 different items in it, then the probability that a specific item is ejected is $1/5$.  

\subsection{The observer}
An \textit{observer} is a block update detector.  This means that an observer can detect when the block it faces experiences a change.  Some examples of how a block might change are crops growing, ice melting, or fire spreading.   These changes occur randomly so a randomizer can be created by detecting these changes with an observer.

\section{Square root of $2$}
The first mathematical constant we will discuss is the square root of two, often denoted as $\sqrt{2}$.  We chose to discuss this number first because it is believed that $\sqrt{2}$ is the first number shown to be irrational.  This means that $\sqrt{2}$ cannot be expressed as the ratio of two integers.
The proof that $\sqrt{2}$ is irrational was first done by the Pythagoreans over 2000 years ago \cite{Kline}.
The proof they gave is an example of a \emph{reductio ad adsurdum}, an argument that establishes a claim by showing that the opposite scenario would lead to contradiction.  
Here is the argument given: let the ratio of hypotenuse to leg of an isosceles right triangle be $a:b$ and let this ratio be expressed in the smallest whole numbers i.e. $a$ and $b$ have no common factors.  By the Pythagorean theorem, we have $a^2 = 2b^2$.  This means that $a^2$ is even and as a result $a$ is even since the square of an odd number is odd.  Since $a$ is even and $a$ and $b$ have no common factors, it follows that $b$ must be odd.  However since $a$ is even, it can be expressed as $a = 2c$ for some integer $c$.  Substituting this into the equation $a^2 = 2b^2$, we have $4c^2 = 2b^2$. It follows that $b$ must also be even by a similar argument as above, but this is a contradiction as $b$ is odd.

Coincidentally, $\sqrt{2}$ is likely one of the first irrational numbers you encountered in your mathematical career.  You probably were first introduced to this number in geometry or trigonometry class when learning about special right triangles.  As a reminder, the ratio of the side lengths of a $45^\circ - 45^\circ - 90^\circ$ right triangle is $1 : 1 : \sqrt{2}$.  In decimal form, the value of $\sqrt{2}$ to four decimal places is $\sqrt{2} = 1.4142$.

\subsection{Approximating square root of $2$ in Minecraft}

\begin{figure}[h]
\begin{center}
\includegraphics[height=5cm,width=5cm]{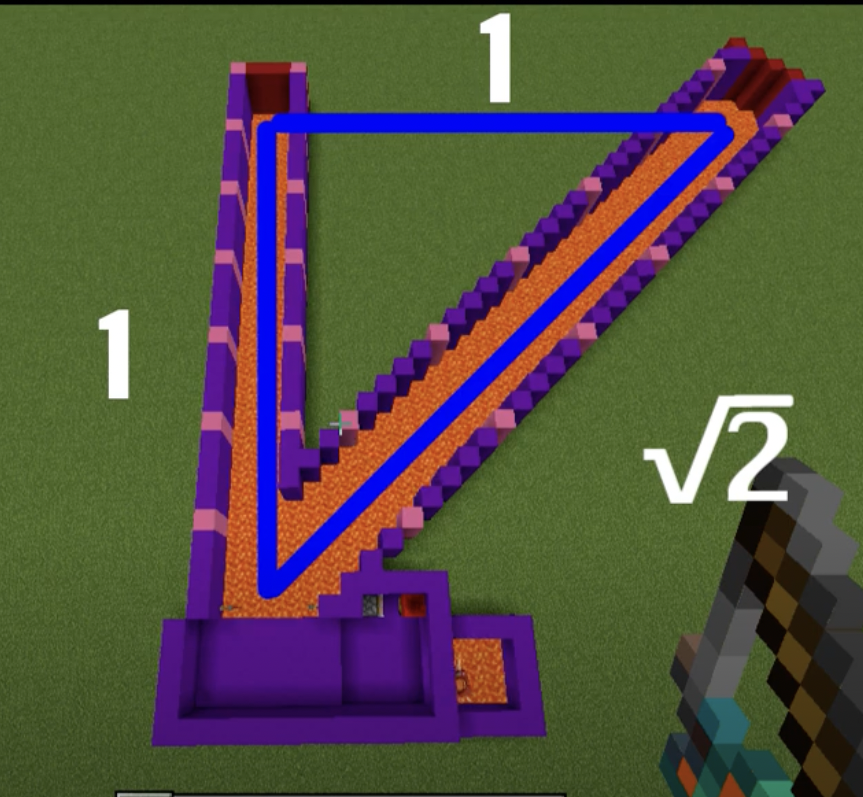}
\caption{The leg and hypotenuse of a large $45^\circ - 45^\circ - 90^\circ$ right triangle.}\label{lavatriangle}
\end{center}
\end{figure}

We will now describe our method for approximating $\sqrt{2}$ in Minecraft.
We started by constructing the leg and hypotenuse of a large $45^\circ - 45^\circ - 90^\circ$ right triangle.  This is one of the easiest shapes to make in Minecraft (see Figure \ref{lavatriangle}) since any block placed in Minecraft must be placed on a grid.
The ratio of the length of the hypotenuse to the length of the leg of the triangle is $\sqrt{2} : 1$. Therefore, to approximate the value of $\sqrt{2}$, we can simply measure the time it takes for a player to travel at a constant rate along one of the legs and the time it takes for that player to travel along the hypotenuse.  Since the length of the hypotenuse is $\sqrt{2}$ times as long as the length of the leg, the ratio of the traveling times should be approximately $\sqrt{2}$.

As mentioned earlier, the timing mechanism can be automated in Minecraft using a hopper timer as hoppers release items at a constant rate.  For both the hypotenuse and the leg, we kept track of how many items the hopper released during the time that the player was traveling. Once these data were collected, taking the ratio of the number of items released while traveling along the hypotenuse to the number of items released while traveling along the leg gives a fractional approximation for $\sqrt{2}$.  

When we completed this experiment, the hopper released 57 items during the time spent traveling along the hypotenuse.  When traveling along the leg of the triangle, the hopper released 41 items.
This gives our approximation as \[\sqrt{2} \approx \frac{57}{41} = 1.3902.\] 
The error for our approximation is $1.70\%.$

\subsubsection{Try this on your own}
There are some easy ways that our approximation can be improved on.  One obvious way to get a more accurate approximation is by making a larger triangle.  If you don't want to spend the time building a larger triangle, another way to improve the approximation is by using a slower method of transportation.  For example, the player can consume a Potion of Slowness before traveling to decrease the player's speed.  Both of these methods will increase the amount of time that the player is traveling.  As a result, the potential error from the hopper timer will be a smaller proportion of the overall time spent traveling.  We encourage you to try this experiment on your own to see how accurate your approximation is.  If you cannot figure out how to make a hopper timer, you can use a stopwatch instead.

The method we used to approximate $\sqrt{2}$ can also be used to approximate the square root of other numbers.  For example, if you want to approximate the value of $\sqrt{5}$, you can make a rectangle whose side lengths are $1$ and $2$ and measure the time it takes to travel along the diagonal and divide it by the time it takes to travel along the shorter leg.  
The reason why $\sqrt{5}$ can be approximated with this method is because $5$ can be expressed as the sum of two perfect squares, i.e. $5 = 1^2 + 2^2$.  However, this is not true of all whole numbers.  For example, this method cannot be used to approximate $\sqrt{7}$ since $7$ cannot be expressed as the sum of two perfect squares.
This leads to the question: what numbers can be expressed as the sum of two squares?  For those of you who are geometry teachers, you can use these experiments as a nice way to introduce basic number theory to geometry students.

\section{$\pi$}

The next mathematical constant we will discuss is $\pi$. The value of $\pi$ to five decimal places is $\pi = 3.14159$. Like $\sqrt{2}$, you very likely encountered $\pi$ in geometry class as it appears in the formulas for the area and circumference of a circle.
For any circle, the ratio of its circumference to its diameter is $\pi$. It was unknown whether or not $\pi$ is an irrational number until 1768 when Johann Heinrich Lambert showed that $\pi$ is indeed irrational \cite{Beckmann}. In fact, $\pi$ is a special type of irrational number known as a \emph{transcendental number}.  This means that $\pi$ is never the root of a non-zero polynomial with integer coefficients. Note that not all irrational numbers are transcendental.  For example, $\sqrt{2}$ is irrational, but it is not transcendental since $\sqrt{2}$ is a root of the polynomial equation $x^2-2=0$.  In the paper where Lambert proved that $\pi$ is irrational, he conjectured that it is also transcendental, but he was unable to prove his conjecture.  It was not until $114$ years later that in 1882 Carl Louis Ferdinand von Lindemann first proved that $\pi$ is transcendental \cite{Beckmann}. 
Despite these properties being unknown until relatively recently, approximations of $\pi$ have existed for millennia.

\begin{figure}[h]
\begin{center}
\begin{tikzpicture}[scale = .75]
\draw[thin] (-4,0) circle (1);
\draw[red, thin] (-4,0) \polygon{3}{1};
\draw[blue, thin ] (-4,0) \polygon{3}{.5};
\draw[black, thin] (0,0) circle (1);
\draw[red, thin] (0,0) \polygon{6}{1};
\draw[blue] (0,0) \polygon{6}{.855};
\draw[black, thin] (4,0) circle (1);
\draw[red, thin] (4,0) \polygon{12}{1};
\draw[blue, thin] (4,0) \polygon{12}{.955};
\end{tikzpicture}
\caption{The circumference of the circle is bounded below by the perimeter of the inscribed polygon (blue) and bounded above by the perimeter of the polygon (red) that circumscribes the circle.}\label{Archimedes Circle}
\end{center}
\end{figure}
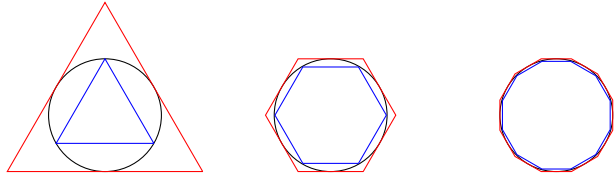

The Greek mathematician Archimedes is credited with the first algorithm to rigorously calculate the value of $\pi$.  Archimedes found lower and upper bounds for the value of $\pi$ by constructing regular polygons inside and outside of the circle (see Figure \ref{Archimedes Circle}). 
He was able to compute the perimeter of these polygons using basic geometry techniques (see \cite{archimedespi} for the specifics).  When using a polygon with $96$ sides, Archimedes found that $3.1408 < \pi < 3.1429$.

\begin{figure}[h]
\begin{center}
\includegraphics[height=5cm,width=5cm]{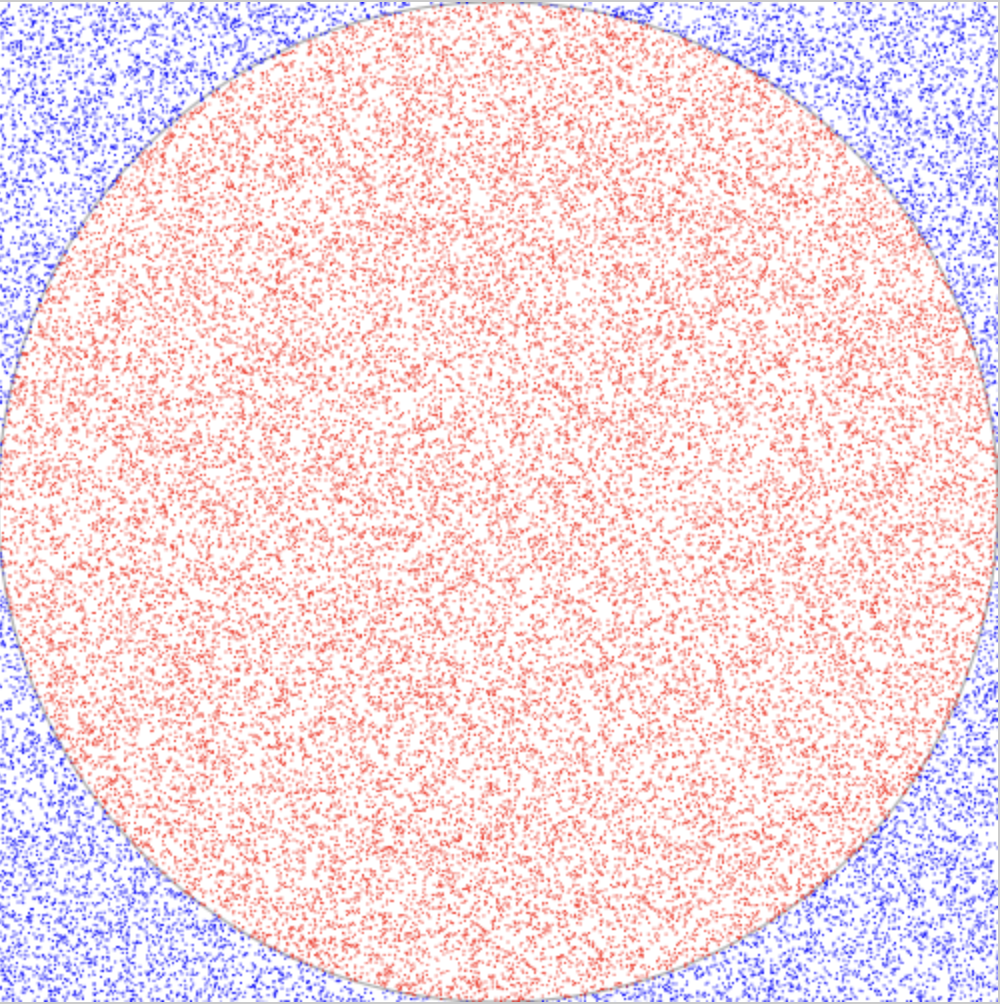}
\caption{After 43270 dots, 33943 in the circle, the approximation of $\pi$ is $\pi \approx 3.13779$, image taken from \cite{Montecarlopic}.}\label{Monte Carlo pi}
\end{center}
\end{figure}

The development of computers led to different methods for computing the digits of $\pi$.  One category of methods is \emph{Monte Carlo methods}. These methods evaluate the results of multiple random trials to approximate values.  One Monte Carlo method, Monte Carlo integration, for approximating the value of $\pi$ is done by drawing a unit circle circumscribed by a square. Then dots are scattered uniformly at random in the square. Since the area of the circle is $\pi$ and the area of the square is $4$, the ratio of dots inside the circle to the total number of dots will approximately equal $\frac{\pi}{4}$ (see Figure \ref{Monte Carlo pi} for an example of Monte Carlo integration).  

\subsection{Approximating $\pi$ in Minecraft}

\begin{figure}[h]
\begin{center}
\includegraphics[height=5cm,width=5cm]{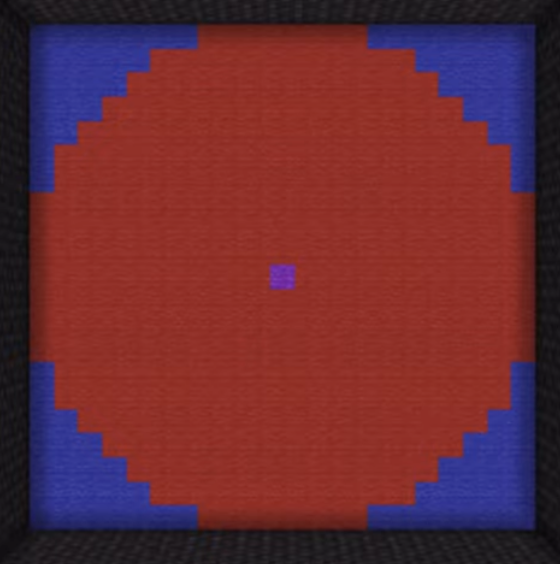}
\caption{An approximation of a circle with radius $11$ in Minecraft}\label{circleMinecraft}
\end{center}
\end{figure}

We recreated Monte Carlo integration to approximate the value of $\pi$ in Minecraft.  Recall that every block in Minecraft is placed on the grid so it is impossible to make a perfect circle.  However, there are many resources online that can approximate the boundary of a circle in Minecraft.  We used the following \href{https://clickspeeder.com/pixel-circle-generator/}{Minecraft circle generator} to make the circle we used (see Figure \ref{circleMinecraft}) \cite{Minecraftcirclegenerator}.

\begin{figure}[h]
\begin{center}
\includegraphics[height=5cm,width=5cm]{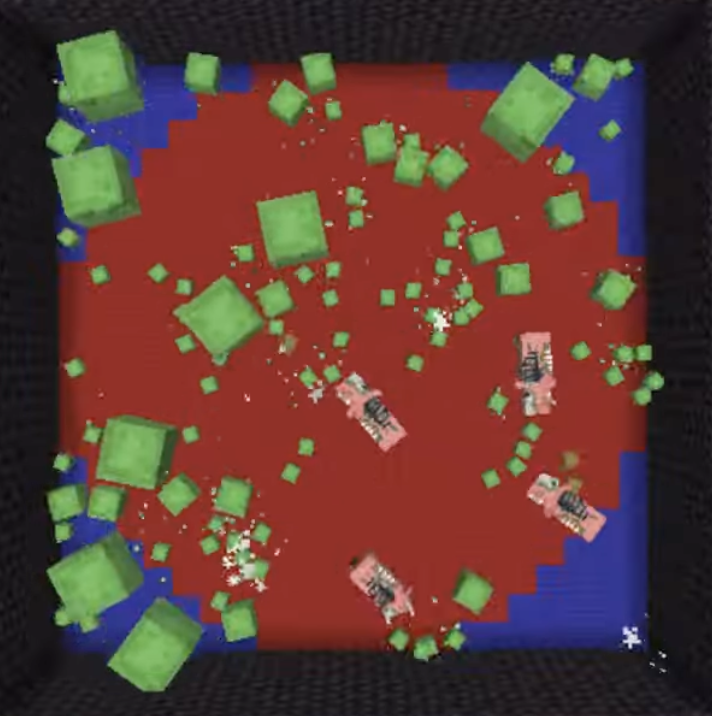}
\caption{Slimes being killed by zoglins.  The location of each death was recorded by hoppers.}\label{Slime Pi}
\end{center}
\end{figure}

The next challenge was to find a way to generate random dots in Minecraft.  To do this we used the behavior of a mob called \emph{the slime}.  We used slimes because unlike other mobs, slimes continue moving when no players are nearby and they change direction at random \cite{slimewiki}.  In fact, most other mobs have a bias to walk in the south-east direction so they would bunch up in the south-east corner of the square \cite{southeastwiki}.  We had a different type of mob (\emph{zoglin}) kill the slimes and used hoppers to keep track of whether or not the slime was killed inside of the circle (see Figure \ref{Slime Pi}). 

In our experiment, a total of $619$ slimes were killed, $508$ of which were killed inside the circle.  This gives our approximation as 
\[\pi \approx 4 \cdot \frac{508}{619} = 3.283.\]
The error for our approximation is $4.49\%.$  The size of our error is not surprising since Monte Carlo methods tend to be slow to converge.  

\subsubsection{Try this on your own}  This method can be improved upon by increasing the size of the circle and increasing the number of slimes that get killed.  Typically in this Monte Carlo method the size of the circle would not affect the accuracy of the approximation, but since you cannot make a perfect circle in Minecraft, increasing the size of the circle will improve the accuracy.  

\begin{figure}[h]
\begin{center}
\includegraphics[height=5cm,width=5cm]{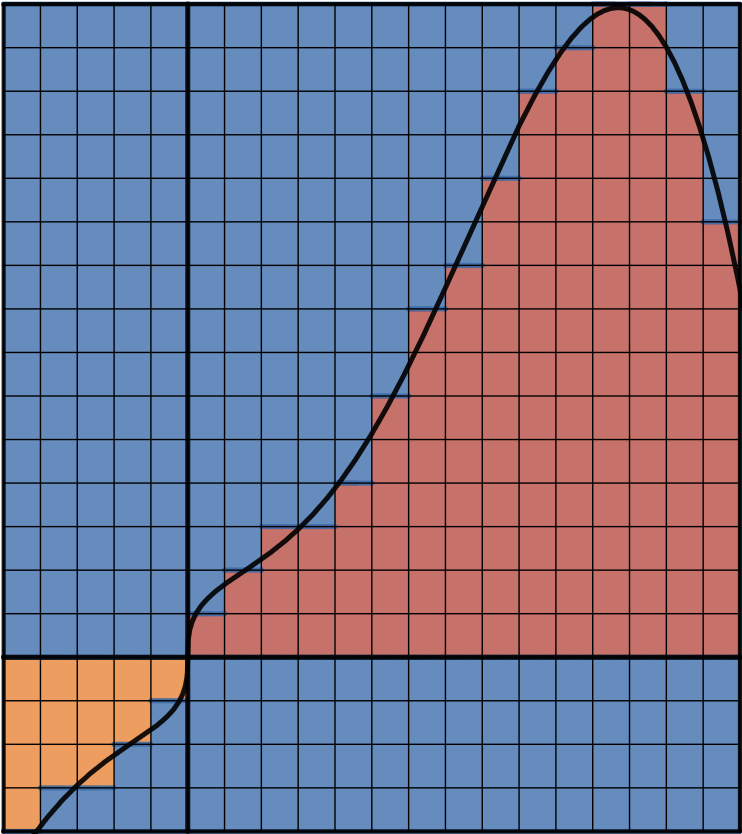}
\caption{The graph of $f(x) = x^2\sin(x) +\sqrt[3]{x}$ and how to make it in Minecraft.}\label{Desmos Graph}
\end{center}
\end{figure}

The techniques we used to approximate the value of $\pi$ can also be used to approximate the values of other definite integrals.  For example, suppose you wanted to use Minecraft to do Monte Carlo integration to approximate the value of $\int_a^bf(x) \ dx$.  With the assistance of \href{https://www.desmos.com/calculator/s80lr06ccg}{this Desmos page} that we created, you can plot out the region between the curve of $y = f(x)$ and the $x$-axis (see Figure \ref{Desmos Graph}) \cite{Desmospage}. 
Recall that the value of a definite integral $\int_a^bf(x) \ dx$ is the net area of the region bounded by the curve $y=f(x)$ and the $x$-axis between $x=a$ and $x=b$.  Therefore, one way to approximate the value of a definite integral in Minecraft is to first find the difference between the number of slimes that die in the region above the $x$-axis with the number of slimes that die in the region below the $x$-axis.  Multiplying this difference by the total area and dividing by the total number of slimes that died will give an approximation for the value of the definite integral.  One function that can be used to approximate the curve $y=f(x)$ in Minecraft is $y = \lfloor f(\lfloor x \rfloor + 0.5) \rceil$ (here $\lfloor x \rceil$ rounds $x$ to its nearest whole number).

This could be a fun experiment for students who are being introduced to integral calculus.  In addition to having students approximate areas using left endpoints, right endpoints, or midpoints, why not have them also approximate areas with Minecraft?

\section{Euler's number}

Another mathematical constant that you have likely come across before is Euler's number $e$.  The value of $e$ to five decimal places is $e = 2.71828$. 
You may recall that $e$ is the base of the natural logarithm and that it is part of the compound interest formulas. It is defined as the following limit:
\begin{equation}\label{e definition}
e = \lim_{n\rightarrow \infty} \left(1+\frac{1}{n}\right)^{n}.
\end{equation}

Although logarithms of base $e$ were being computed as early 1618, the use of the symbol $e$ was not referenced in these calculations. The first so-called discovery of $e$ came about through work of Jacob Bernoulli in 1638 while studying continuous compound interest. In doing so, he tried to evaluate the limit seen in Equation \ref{e definition}. Using the binomial theorem, he was able to show that $2 < e < 3$, but at this time $e$ still did not have a name nor a closer approximation \cite{thenumbere}.

It was Euler who ended up making the connection between logarithms and the number $e$. He  evaluated the limit seen in Equation \ref{e definition} and gave its value the symbol $e$.   In 1737 he proved that $e$ is  irrational.  In 1873, Charles Hermite showed that $e$ is transcendental.  In fact, $e$ is the first number to be proven transcendental without having been specifically constructed to be transcendental \cite{Hermitebio}.

For our purposes it is necessary to use the expression of $e$ given by Euler in 1748. He published \textit{Introductio in Analysin infinitorum} \cite{Euler}, where he showed that $e$ can be expressed as \[e = 1 + \frac{1}{1!} + \frac{1}{2!}+ \frac{1}{3!} + \cdots .\]

We will now consider the function $f(x) = e^x$.  This function can be expressed as its Maclaurin series 
\[e^{x} = \sum_{k=0}^{\infty} \frac{x^{k}}{k!}.\] Note that evaluating at $x=-1$ gives the following alternating series expansion for $\frac{1}{e}$: 
\begin{equation}\label{1/e}
\frac{1}{e} = \sum_{k=0}^{\infty} \frac{(-1)^{k}}{k!}.
\end{equation}
We will see that the $n$th partial sum of this expression is the solution to a specific counting problem.  We will now describe this problem.

Let $[n] = \{1,2,3,...,n\}$.

\begin{definition} 
A \emph{permutation of $[n]$} is an arrangement of $[n]$ in a definite order. 
\end{definition}

The permutations of $[n]$ can be thought of as the linear orderings of the numbers $1$ through $n$.  For example, the permutations of $[3]$ are $123, 132, 213, 231, 312,$ and $ 321$.  The total number of permutations of $[n]$ is $n \cdot (n-1) \cdots 2 \cdot 1$.  This product is traditionally denoted by $n!$.  

\begin{definition}
A \textit{derangement} is a permutation of $[n]$ that has no fixed points.
\end{definition}

In other words, if $\omega$ is a permutation of $[n]$, then $\omega$ is a derangement if $\omega_i \neq i$ for all $i \in [n]$.  For example, consider the following permutation of $[6]$: $\omega = 324165$.  
This is not a derangement because the number 2 is in the second position, i.e. $\omega_{2} = 2$. However, the permutation $\nu=431562$ is a derangement since $\nu_{i} \neq i$ for all $i \in [6]$. For more on permutations and derangements, see \cite{Stanley}.

We will denote the number of derangements of $[n]$ by $D(n)$. It can be shown that 
\begin{equation}\label{D(n)}
D(n) = \sum_{k=0}^{n} \binom{n}{k}(-1)^{k}(n-k)! = \sum_{k=0}^{n} (-1)^{k}\frac{n!(n-k)!}{k!(n-k)!} = \sum_{k=0}^{n} \frac{(-1)^{k}n!}{k!}.
\end{equation}

Comparing Equations \ref{1/e} and \ref{D(n)}, you can see that $\frac{D(n)}{n!}$ gives the $n$th partial sum of Equation \ref{1/e}.  
As a result, we can see that $\frac{1}{e}$ is approximately equal to the probability that a random permutation is a derangement.  In particular, $e \approx \frac{n!}{D(n)}$.

\subsection{Approximating Euler's number in Minecraft}

To approximate $e$ in Minecraft, we made a machine that:
\begin{enumerate}
\item produces a permutation and 
\item checks whether or not it is a derangement.
\end{enumerate}
Once the machine was created, we then had the machine run many times to generate a large enough sample.  We will now describe how we made said machine.

As mentioned above, a dropper can be used as a randomizer.  
Since a dropper can hold up to 9 different items, we can use its random ejection mechanism to create permutations of $[9]$.  Each block in the dropper would correspond to a different number in $[9]$. The order that the blocks are ejected can be read as a permutation.

There are ways in Minecraft to automatically check which block was ejected by the dropper.  The details of how to do this are beyond the scope of this article (for those who are interested see \cite{impulseSV}).
As a result, one can make a machine that checks whether or not the permutation that was created is a derangement.  This is done by checking if the block corresponding to a certain number was ejected in that position.  If each of the $9$ blocks were ejected in a position that didn't correspond to their number, then the permutation is a derangement.

The proportion of the permutations that we generated that were derangements is approximately $\frac{1}{e}$.
That is to say 
\[e \approx \frac{\text{number of permutations}}{\text{number of derangements}}.\] We generated 647 permutations and 238 of these were derangements. This gives us our approximation of $e$ as \[e \approx \frac{647}{238} = 2.71849.\] The error for our approximation was about $0.00766 \%.$  We did get lucky with the accuracy based on the number of permutations we generated.  However, if we let the machine run indefinitely, the error for the approximation of $\frac{1}{e}$ would be less than $\frac{1}{10!}$.

\subsubsection{Try this on your own} A permutation $\omega$ is said to be \emph{alternating} if its entries alternately rise and descend. That is $\omega_1<\omega_2>\omega_3<...$.  
For example, the permutation $1423$ is alternating, but the permutation $1342$ is not since $\omega_2 \ngtr \omega_3$.
If we let $A_n$ denote the number of alternating permutation of $[n]$, then Andr\'{e}'s Theorem \cite{Andrestheorem} states that 
\[\sum_{n=0}^\infty\frac{A_n}{n!} = \sec(1)+\tan(1).\]
We're not sure why you'd want to do this, but this means you can approximate the value of $\sec(1)+\tan(1)$ using Minecraft.  Like we did for the approximation of $e$, you can use droppers to generate permutations.  Instead of having one dropper make a permutation of $[9]$, you can use $9$ (or fewer) droppers to generate permutations of different sizes. For each of the permutations, you can make a machine that checks if it's alternating.  You would then compute the proportion of alternating permutations for each size and add these proportions together.  This gives an approximation of the $9$th partial sum of the series.  Please contact the authors if you complete this experiment and let us know your results.

\section{Ap\'{e}ry's constant}

The final mathematical constant that we will discuss is one that you likely have not seen before and if you have encountered it, you may not know that it has a name. Ap\'{e}ry's constant, denoted $\zeta(3),$ is defined to be the sum of the reciprocals of the positive cubes. That is 
\begin{equation}
\zeta(3) = \sum_{n=1}^{\infty} \frac{1}{n^{3}} = \lim_{n\rightarrow \infty} \left( \frac{1}{1^{3}} + \frac{1}{2^{3}} + \cdots + \frac{1}{n^{3}}\right). \nonumber
\end{equation}
The reason why Ap\'{e}ry's constant is denoted by $\zeta(3)$ is that  it is the value of the Riemann zeta function evaluated at $s=3$.
In general, the Riemann zeta functions is defined as
\[
\zeta(s) = \sum_{n=1}^{\infty} \frac{1}{n^{s}}.
\]

Euler proved the following product formula for the  Riemann zeta function:

\begin{equation}
\zeta(s) = 
\prod_{p= \text{prime}} \Bigl(1 + \frac{1}{p^s} + \Bigl(\frac{1}{p^s}\Bigr)^2 + \dots \Bigr). \label{Zeta}
\end{equation}


The value of Ap\'{e}ry's constant to $5$ decimal places is $\zeta(3) = 1.20205.$ It was not until 1979 that Roger Ap\'{e}ry (from whom the number gets its name) showed $\zeta(3)$ is irrational \cite{aperyirrational}. It is still an open question whether or not the number is transcendental. There are various representation of Ap\'{e}ry's constant as series and as integrals. Some of these representations are much more complicated than others. However, 
we will make use of the fact that value of Ap\'{e}ry's constant can be determined probabilistically.  The reciprocal of Ap\'{e}ry's constant is the probability that any three positive integers chosen uniformly at random will be relatively prime. We will now give a brief sketch of why this is the case.

For three positive integers to be relatively prime, there is no prime number that divides evenly into all three of the numbers.  For example, $6, 9,$ and $ 21$ are not relatively prime since they are each divisible by $3$.  For a prime $p$, the probability that $p$ divides a random whole number is $\frac{1}{p}$.  Therefore $p$ divides all three numbers with probability $\frac{1}{p^3}$.  This means that the probability that at least one of the three numbers is not divisible by $p$ is $(1 - \frac{1}{p^3})$. Let $P_3$ denote the probability that the three randomly selected positive integers are relatively prime.  It follows that  

\begin{align}
P_3 & = \prod_{p= \text{prime}} \Bigl(1 - \frac{1}{p^3}\Bigr) \nonumber \\
 &= \prod_{p= \text{prime}} \Biggl(\frac{1}{1 - \frac{1}{p^3}}\Biggr)^{-1} \nonumber \\
 & = \prod_{p= \text{prime}} \Bigl(1 + \frac{1}{p^3} + \Bigl(\frac{1}{p^3}\Bigr)^2 + \dots \Bigr)^{-1}. \label{Probability}
\end{align}
Comparing Equations \ref{Zeta} and \ref{Probability}, we see that $P_3 = \zeta(3)^{-1}$.

\subsection{Approximating Ap\'{e}ry's constant in Minecraft}

In order to approximate Ap\'{e}ry's constant in Minecraft, we repeatedly generated sets of three random numbers, called triplets, and manually checked if these numbers are relatively prime. 
The proportion of triplets that are relatively prime will be approximately $\zeta(3)^{-1}$.

As mentioned in the Minecraft Mechanics section above, an observer is a block in Minecraft that can detect when the block it is facing experiences a change of state. Many blocks in Minecraft change their state at random intervals.  In general, every  $0.05$ seconds, the game randomly selects 3 blocks (we increased this value to speed up the experiment) in a $16 \times 16 \times 16$ cube to change the state of.  If the selected blocks have the ability to change states, they will change their state with some predetermined probability.  For more details on how blocks change states in Minecraft, see \cite{tickwiki}.

\begin{figure}[h]
\begin{center}
\includegraphics[height=5cm,width=5cm]{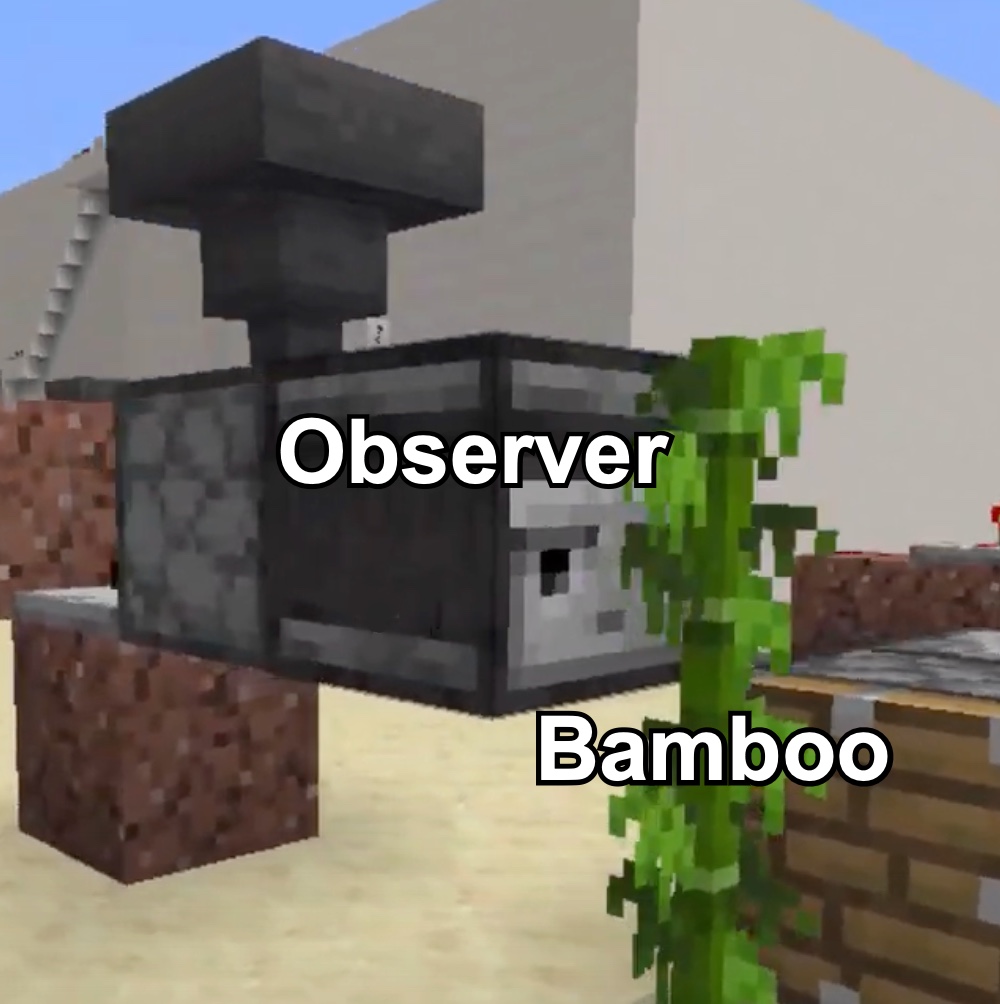}
\caption{The observer detects the exact moment that the bamboo grows.}\label{Observer}
\end{center}
\end{figure}

To generate a triplet of random numbers, we had three observers each facing their own bamboo plant.  We used a hopper timer to record how long it took each bamboo plant to change states (see Figure \ref{Observer}).  It should be noted that the random numbers that were generated did not follow a uniform distribution, instead they followed the negative binomial distribution.

In our experiment, we collected $70$ triplets of random numbers.   We found that $58$ of the triplets were relatively prime.  This gives us our approximation of $\zeta(3)$ as \[\zeta(3) \approx \frac{70}{58} = 1.2069.\]
The error for our approximation is about $0.4 \%.$  The numbers we generated varied between a minimum of $3$ and a maximum of $838$ so our method did a better than we expected job at getting a wide variety of numbers.

\subsubsection{Try this on your own} 
Recall that the probability that three positive integers chosen uniformly at random will be relatively prime is $\zeta(3)^{-1}$.  In general $P_m$, the probability that $m$ positive integers chosen uniformly at random will be relatively prime, is $\zeta(m)^{-1}$.  This means you can use the techniques described above to approximate the values of $\zeta(m)$ for various values of $m$.  In particular, you can approximate the value of $\pi$ raised to any even power since $\zeta(2k)$ is always a rational multiple of $\pi^{2k}$.   For example, since $\zeta(2) = \frac{\pi^2}{6}$, it follows that $\pi^2 = 6 P_2^{-1}$.  This means that you can approximate the value of $\pi^2$ by generating pairs of numbers and checking whether or not they are relatively prime.  For those looking for a challenge, try to automate the checking process by making a machine in Minecraft that finds the greatest common divisor of the two positive integers.  If the greatest common divisor is $1$, then the numbers are relatively prime.  Making this machine will likely be difficult, but not out of the range of what is possible to do in Minecraft.
 
\section{Conclusion}

We hope this paper gives a small sample of what is possible in the world of math and Minecraft and serves as an inspiration to explore complex mathematical topics in fun and interesting ways. While we chose to use Minecraft to approximate irrational numbers, we believe there are many settings where this type of experimentation is possible. We would love to hear from anyone who finds inspiration from this paper and does some exploring of their own!



\end{document}